# Marshall's lemma for convex density estimation

**Lutz Dümbgen**[1], **Kaspar Rufibach**[1] **and Jon A. Wellner**[2]

*University of Bern and University of Washington*

**Abstract:** Marshall's [*Nonparametric Techniques in Statistical Inference* (1970) 174–176] lemma is an analytical result which implies $\sqrt{n}$–consistency of the distribution function corresponding to the Grenander [*Skand. Aktuarietidskr.* **39** (1956) 125–153] estimator of a non-decreasing probability density. The present paper derives analogous results for the setting of convex densities on $[0,\infty)$.

## 1. Introduction

Let $\mathbb{F}$ be the empirical distribution function of independent random variables $X_1, X_2, \ldots, X_n$ with distribution function $F$ and density $f$ on the halfline $[0,\infty)$. Various shape restrictions on $f$ enable consistent nonparametric estimation of it without any tuning parameters (e.g. bandwidths for kernel estimators).

The oldest and most famous example is the Grenander estimator $\hat{f}$ of $f$ under the assumption that $f$ is non-increasing. Denoting the family of all such densities by $\mathcal{F}$, the Grenander estimator may be viewed as the maximum likelihood estimator,

$$\hat{f} = \operatorname{argmax}\Big\{ \int \log h \, d\mathbb{F} : h \in \mathcal{F} \Big\},$$

or as a least squares estimator,

$$\hat{f} = \operatorname{argmin}\Big\{ \int_0^\infty h(x)^2 dx - 2 \int h \, d\mathbb{F} : h \in \mathcal{F} \Big\};$$

cf. Robertson et al. [5]. Note that if $\mathbb{F}$ had a square-integrable density $\mathbb{F}'$, then the preceding argmin would be identical with the minimizer of $\int_0^\infty (h - \mathbb{F}')(x)^2 \, dx$ over all non-increasing probability densities $h$ on $[0,\infty)$.

A nice property of $\hat{f}$ is that the corresponding distribution function $\hat{F}$,

$$\hat{F}(r) := \int_0^r \hat{f}(x) \, dx,$$

is automatically $\sqrt{n}$–consistent. More precisely, since $\hat{F}$ is the least concave majorant of $\mathbb{F}$, it follows from Marshall's [4] lemma that

$$\|\hat{F} - F\|_\infty \ \leq \ \|\mathbb{F} - F\|_\infty.$$

A more refined asymptotic analysis of $\hat{F} - \mathbb{F}$ has been provided by Kiefer and Wolfowitz [3].

---

[1]Institute of Math. Statistics and Actuarial Science, University of Bern, Switzerland, e-mail: lutz.duembgen@stat.unibe.ch; kaspar.rufibach@stat.unibe.ch
[2]Dept. of Statistics, University of Washington, Seattle, USA, e-mail: jaw@stat.washington.edu

## 2. Convex densities

Now we switch to the estimation of a convex probability density $f$ on $[0, \infty)$. As pointed out by Groeneboom et al. [2], the nonparametric maximum likelihood estimator $\hat{f}_{ml}$ and the least squares estimator $\hat{f}_{ls}$ are both well-defined and unique, but they are not identical in general. Let $\mathcal{K}$ denote the convex cone of all convex and integrable functions $g$ on $[0, \infty)$. (All functions within $\mathcal{K}$ are necessarily nonnegative and non-increasing.) Then

$$\hat{f}_{ml} = \underset{h \in \mathcal{K}}{\operatorname{argmax}} \Big( \int \log h \, d\mathbb{F} - \int_0^\infty h(x) \, dx \Big),$$

$$\hat{f}_{ls} = \underset{h \in \mathcal{K}}{\operatorname{argmin}} \Big( \int_0^\infty h(x)^2 dx - 2 \int h \, d\mathbb{F} \Big).$$

Both estimators have the following property:

**Proposition 1.** *Let $\hat{f}$ be either $\hat{f}_{ml}$ or $\hat{f}_{ls}$. Then $\hat{f}$ is piecewise linear with*

- *at most one knot in each of the intervals $(X_{(i)}, X_{(i+1)})$, $1 \le i < n$,*
- *no knot at any observation $X_i$, and*
- *precisely one knot within $(X_{(n)}, \infty)$.*

The estimators $\hat{f}_{ml}$, $\hat{f}_{ls}$ and their distribution functions $\hat{F}_{ml}$, $\hat{F}_{ls}$ are completely characterized by Proposition 1 and the next proposition.

**Proposition 2.** *Let $\Delta$ be any function on $[0, \infty)$ such that $\hat{f}_{ml} + t\Delta \in \mathcal{K}$ for some $t > 0$. Then*

$$\int \frac{\Delta}{\hat{f}_{ml}} \, d\mathbb{F} \ \le \ \int \Delta(x) \, dx.$$

*Similarly, let $\Delta$ be any function on $[0, \infty)$ such that $\hat{f}_{ls} + t\Delta \in \mathcal{K}$ for some $t > 0$. Then*

$$\int \Delta \, d\mathbb{F} \ \le \ \int \Delta \, d\hat{F}_{ls}.$$

In what follows we derive two inequalities relating $\hat{F} - F$ and $\mathbb{F} - F$, where $\hat{F}$ stands for $\hat{F}_{ml}$ or $\hat{F}_{ls}$:

**Theorem 1.**

(1) $$\inf_{[0,\infty)} (\hat{F}_{ml} - F) \ge \frac{3}{2} \inf_{[0,\infty)} (\mathbb{F} - F) - \frac{1}{2} \sup_{[0,\infty)} (\mathbb{F} - F),$$

(2) $$\big\| \hat{F}_{ls} - F \big\|_\infty \le 2 \, \big\| \mathbb{F} - F \big\|_\infty.$$

Both results rely on the following lemma:

**Lemma 1.** *Let $F, \hat{F}$ be continuous functions on a compact interval $[a, b]$, and let $\mathbb{F}$ be a bounded, measurable function on $[a, b]$. Suppose that the following additional assumptions are satisfied:*

(3) $\quad\quad\quad\quad \hat{F}(a) = \mathbb{F}(a) \ \text{and} \ \hat{F}(b) = \mathbb{F}(b),$

(4) $\quad\quad\quad\quad \hat{F} \ \text{has a linear derivative on} \ (a, b),$

(5) $\quad\quad\quad\quad F \ \text{has a convex derivative on} \ (a, b),$

(6) $\quad\quad\quad\quad \int_r^b \hat{F}(y) \, dy \ \le \ \int_r^b \mathbb{F}(y) \, dy \quad \text{for all} \ r \in [a, b].$

*Then*

$$\sup_{[a,b]} (\hat{F} - F) \leq \frac{3}{2} \sup_{[a,b]} (\mathbb{F} - F) - \frac{1}{2}(\mathbb{F} - F)(b).$$

*If condition (6) is replaced with*

(7) $$\int_a^r \hat{F}(x)\,dx \geq \int_a^r \mathbb{F}(x)\,dx \quad \text{for all } r \in [a,b],$$

*then*

$$\inf_{[a,b]} (\hat{F} - F) \geq \frac{3}{2} \inf_{[a,b]} (\mathbb{F} - F) - \frac{1}{2}(\mathbb{F} - F)(a).$$

The constants $3/2$ and $1/2$ are sharp. For let $[a,b] = [0,1]$ and define

$$F(x) := \begin{cases} x^2 - c & \text{for } x \geq \epsilon, \\ (x/\epsilon)(\epsilon^2 - c) & \text{for } x \leq \epsilon, \end{cases}$$
$$\hat{F}(x) := 0,$$
$$\mathbb{F}(x) := 1\{0 < x < 1\}(x^2 - 1/3)$$

for some constant $c \geq 1$ and some small number $\epsilon \in (0, 1/2]$. One easily verifies conditions (3)–(6). Moreover,

$$\sup_{[0,1]} (\hat{F} - F) = c - \epsilon^2, \quad \sup_{[0,1]} (\mathbb{F} - F) = c - 1/3 \quad \text{and} \quad (\mathbb{F} - F)(1) = c - 1.$$

Hence the upper bound $(3/2)\sup(\mathbb{F} - F) - (1/2)(\mathbb{F} - F)(1)$ equals $\sup(\hat{F} - F) + \epsilon^2$ for any $c \geq 1$. Note the discontinuity of $\mathbb{F}$ at $0$ and $1$. However, by suitable approximation of $\mathbb{F}$ with continuous functions one can easily show that the constants remain optimal even under the additional constraint of $\mathbb{F}$ being continuous.

*Proof of Lemma 1.* We define $G := \hat{F} - F$ with derivative $g := G'$ on $(a,b)$. It follows from (3) that

$$\max_{\{a,b\}} G = \max_{\{a,b\}} (\mathbb{F} - F) \leq \frac{3}{2} \sup_{[a,b]} (\mathbb{F} - F) - \frac{1}{2}(\mathbb{F} - F)(b).$$

Therefore it suffices to consider the case that $G$ attains its maximum at some point $r \in (a,b)$. In particular, $g(r) = 0$. We introduce an auxiliary linear function $\bar{g}$ on $[r,b]$ such that $\bar{g}(r) = 0$ and

$$\int_r^b \bar{g}(y)\,dy = \int_r^b g(y)\,dy = G(b) - G(r).$$

Note that $g$ is concave on $(a,b)$ by (4)–(5). Hence there exists a number $y_o \in (r,b)$ such that

$$g - \bar{g} \begin{cases} \geq 0 \text{ on } [r, y_o], \\ \leq 0 \text{ on } [y_o, b). \end{cases}$$

This entails that

$$\int_r^y (g - \bar{g})(u)\,du = -\int_y^b (g - \bar{g})(u)\,du \geq 0 \quad \text{for any } y \in [r,b].$$

Consequently,

$$G(y) = G(r) + \int_r^y g(u)\, du$$
$$\geq G(r) + \int_r^y \bar{g}(u)\, du$$
$$= G(r) + \frac{(y-r)^2}{(b-r)^2}[G(b) - G(r)],$$

so that

$$\int_r^b G(y)\, dy \geq (b-r)G(r) + \frac{G(b) - G(r)}{(b-r)^2} \int_r^b (y-r)^2\, dy$$
$$= (b-r)\left[\frac{2}{3}G(r) + \frac{1}{3}G(b)\right]$$
$$= (b-r)\left[\frac{2}{3}G(r) + \frac{1}{3}(\mathbb{F} - F)(b)\right].$$

On the other hand, by assumption (6),

$$\int_r^b G(y)\, dy \leq \int_r^b (\mathbb{F} - F)(y)\, dy \leq (b-r)\sup_{[a,b]}(\mathbb{F} - F).$$

This entails that

$$G(r) \leq \frac{3}{2}\sup_{[a,b]}(\mathbb{F} - F) - \frac{1}{2}(\mathbb{F} - F)(b).$$

If (6) is replaced with (7), then note first that

$$\min_{\{a,b\}} G = \min_{\{a,b\}}(\mathbb{F} - F) \geq \frac{3}{2}\min_{\{a,b\}}(\mathbb{F} - F) - \frac{1}{2}(\mathbb{F} - F)(a).$$

Therefore it suffices to consider the case that $G$ attains its minimum at some point $r \in (a,b)$. Now we consider a linear function $\bar{g}$ on $[a,r]$ such that $\bar{g}(r) = 0$ and

$$\int_a^r \bar{g}(x)\, dx = \int_a^r g(x)\, dx = G(r) - G(a).$$

Here concavity of $g$ on $(a,b)$ entails that

$$\int_a^x (g - \bar{g})(u)\, du = -\int_x^r (g - \bar{g})(u)\, du \leq 0 \quad \text{for any } x \in [a,r],$$

so that

$$G(x) = G(r) - \int_x^r g(u)\, du$$
$$\leq G(r) - \int_x^r \bar{g}(u)\, du$$
$$= G(r) - \frac{(r-x)^2}{(r-a)^2}[G(r) - G(a)].$$

Consequently,

$$\int_a^r G(x)\, dx \leq (r-a)G(r) - \frac{G(r) - G(a)}{(r-a)^2}\int_a^r (r-x)^2\, dx$$
$$= (r-a)\left[\frac{2}{3}G(r) + \frac{1}{3}(\mathbb{F} - F)(a)\right],$$

whereas
$$\int_a^r G(x)\,dx \;\geq\; \int_a^r (\mathbb{F} - F)(x)\,dx \;\geq\; (r-a)\inf_{[a,b]}(\mathbb{F} - F),$$
by assumption (7). This leads to
$$G(r) \;\geq\; \frac{3}{2}\inf_{[a,b]}(\mathbb{F}-F) - \frac{1}{2}(\mathbb{F}-F)(a). \qquad \square$$

$\square$

*Proof of Theorem 1.* Let $0 =: t_0 < t_1 < \cdots < t_m$ be the knots of $\hat{f}$, including the origin. In what follows we derive conditions (3)–(5) and (6/7) of Lemma 1 for any interval $[a,b] = [t_k, t_{k+1}]$ with $0 \leq k < m$. For the reader's convenience we rely entirely on Proposition 2. In case of the least squares estimator, similar inequalities and arguments may be found in Groeneboom et al. [2].

Let $0 < \epsilon < \min_{1 \leq i \leq m}(t_i - t_{i-1})/2$. For a fixed $k \in \{1, \ldots, m\}$ we define $\Delta_1$ to be continuous and piecewise linear with knots $t_{k-1} - \epsilon$ (if $k > 1$), $t_{k-1}$, $t_k$ and $t_k + \epsilon$. Namely, let $\Delta_1(x) = 0$ for $x \notin (t_{k-1} - \epsilon, t_k + \epsilon)$ and

$$\Delta_1(x) := \begin{cases} \hat{f}_{ml}(x) & \text{if } \hat{f} = \hat{f}_{ml} \\ 1 & \text{if } \hat{f} = \hat{f}_{ls} \end{cases} \quad \text{for } x \in [t_{k-1}, t_k].$$

This function $\Delta_1$ satisfies the requirements of Proposition 2. Letting $\epsilon \searrow 0$, the function $\Delta_1(x)$ converges pointwise to

$$\begin{cases} 1\{t_{k-1} \leq x \leq t_k\}\hat{f}_{ml}(x) & \text{if } \hat{f} = \hat{f}_{ml}, \\ 1\{t_{k-1} \leq x \leq t_k\} & \text{if } \hat{f} = \hat{f}_{ls}, \end{cases}$$

and the latter proposition yields the inequality

$$\mathbb{F}(t_k) - \mathbb{F}(t_{k-1}) \;\leq\; \hat{F}(t_k) - \hat{F}(t_{k-1}).$$

Similarly let $\Delta_2$ be continuous and piecewise linear with knots at $t_{k-1}$, $t_{k-1} + \epsilon$, $t_k - \epsilon$ and $t_k$. Precisely, let $\Delta_2(x) := 0$ for $x \notin (t_{k-1}, t_k)$ and

$$\Delta_2(x) := \begin{cases} -\hat{f}_{ml}(x) & \text{if } \hat{f} = \hat{f}_{ml} \\ -1 & \text{if } \hat{f} = \hat{f}_{ls} \end{cases} \quad \text{for } x \in [t_{k-1} + \epsilon, t_k - \epsilon].$$

The limit of $\Delta_2(x)$ as $\epsilon \searrow 0$ equals

$$\begin{cases} -1\{t_{k-1} < x < t_k\}\hat{f}_{ml}(x) & \text{if } \hat{f} = \hat{f}_{ml}, \\ -1\{t_{k-1} < x < t_k\} & \text{if } \hat{f} = \hat{f}_{ls}, \end{cases}$$

and it follows from Proposition 2 that

$$\mathbb{F}(t_k) - \mathbb{F}(t_{k-1}) \;\geq\; \hat{F}(t_k) - \hat{F}(t_{k-1}).$$

This shows that $\mathbb{F}(t_k) - \mathbb{F}(t_{k-1}) = \hat{F}(t_k) - \hat{F}(t_{k-1})$ for $k = 1, \ldots, m$. Since $\hat{F}(0) = 0$, one can rewrite this as

(8) $$\mathbb{F}(t_k) \;=\; \hat{F}(t_k) \quad \text{for } k = 0, 1, \ldots, m.$$

Now we consider first the maximum likelihood estimator $\hat{f}_{ml}$. For $0 \leq k < m$ and $r \in (t_k, t_{k+1}]$ let $\Delta(x) := 0$ for $x \notin (t_k - \epsilon, r)$, let $\Delta$ be linear on $[t_k - \epsilon, t_k]$, and let $\Delta(x) := (r - x)\hat{f}_{ml}(x)$ for $x \in [t_k, r]$. One easily verifies, that this function $\Delta$ satisfies the conditions of Proposition 2, too, and with $\epsilon \searrow 0$ we obtain the inequality

$$\int_{t_k}^{r} (r - x)\,\mathbb{F}(dx) \leq \int_{t_k}^{r} (r - x)\,\hat{F}(dx).$$

Integration by parts (or Fubini's theorem) shows that the latter inequality is equivalent to

$$\int_{t_k}^{r} (\mathbb{F}(x) - \mathbb{F}(t_k))\,dx \leq \int_{t_k}^{r} (\hat{F}(x) - \hat{F}(t_k))\,dx.$$

Since $\mathbb{F}(t_k) = \hat{F}(t_k)$, we end up with

$$\int_{t_k}^{r} \mathbb{F}(x)\,dx \leq \int_{t_k}^{r} \hat{F}(x)\,dx \quad \text{for } k = 0, 1, \ldots, m-1 \text{ and } r \in (t_k, t_{k+1}].$$

Hence we may apply Lemma 1 and obtain (1).

Finally, let us consider the least squares estimator $\hat{f}_{ls}$. For $0 \leq k < m$ and $r \in (t_k, t_{k+1}]$ let $\Delta(x) := 0$ for $x \notin (t_k - \epsilon, r)$, let $\Delta$ be linear on $[t_k - \epsilon, t_k]$ as well as on $[t_k, r]$ with $\Delta(t_k) := r - t_k$. Then applying Proposition 2 and letting $\epsilon \searrow 0$ yields

$$\int_{t_k}^{r} (r - x)\,\mathbb{F}(dx) \leq \int_{t_k}^{r} (r - x)\,\hat{F}(dx),$$

so that

$$\int_{t_k}^{r} \mathbb{F}(x)\,dx \leq \int_{t_k}^{r} \hat{F}(x)\,dx \quad \text{for } k = 0, 1, \ldots, m-1 \text{ and } r \in (t_k, t_{k+1}].$$

Thus it follows from Lemma 1 that

$$\inf_{[0,\infty)} (\hat{F} - F) \geq \frac{3}{2} \inf_{[0,\infty)} (\mathbb{F} - F) - \frac{1}{2} \sup_{[0,\infty)} (\mathbb{F} - F) \geq -2\,\|\mathbb{F} - F\|_{\infty}.$$

Alternatively, for $1 \leq k \leq m$ and $r \in [t_{k-1}, t_k)$ let $\Delta(x) := 0$ for $x \notin (r, t_k + \epsilon)$, let $\Delta$ be linear on $[r, t_k]$ as well as on $[t_k, t_k + \epsilon]$ with $\Delta(t_k) := -(t_k - r)$. Then applying Proposition 2 and letting $\epsilon \searrow 0$ yields

$$\int_{r}^{t_k} (t_k - x)\,\mathbb{F}(dx) \geq \int_{r}^{t_k} (t_k - x)\,\hat{F}(dx),$$

so that

$$\int_{r}^{t_k} \mathbb{F}(x)\,dx \geq \int_{t_k}^{r} \hat{F}(x)\,dx \quad \text{for } k = 1, 2, \ldots, m \text{ and } r \in [t_{k-1}, t_k).$$

Hence it follows from Lemma 1 that

$$\sup_{[0,\infty)} (\hat{F} - F) \leq \frac{3}{2} \sup_{[0,\infty)} (\mathbb{F} - F) - \frac{1}{2} \inf_{[0,\infty)} (\mathbb{F} - F) \leq 2\,\|\mathbb{F} - F\|_{\infty}. \qquad \square$$

**Acknowledgment.** The authors are grateful to Geurt Jongbloed for constructive comments and careful reading.